\newcount\ite\ite=1\def\0{\global\ite=1\1}
\def\1{\item{\rm(\romannumeral\the\ite)}\advance\ite1\quad}
\def\phi{\varphi}

\documentclass[12pt,twoside]{article}
\usepackage{amssymb}
\usepackage{setspace}

\font\teneufm=eufm10 scaled \magstep1
\font\seveneufm=eufm7 scaled \magstep1
\font\fiveeufm=eufm5  scaled \magstep1
\newfam\eufmfam

\textfont\eufmfam=\teneufm
\scriptfont\eufmfam=\seveneufm
\scriptscriptfont\eufmfam=\fiveeufm

\newfam\msbfam
\font\tenmsb=msbm10 scaled \magstep1  \textfont\msbfam=\tenmsb
\font\sevenmsb=msbm7 scaled \magstep1 \scriptfont\msbfam=\sevenmsb
\font\fivemsb=msbm5 scaled \magstep1  \scriptscriptfont\msbfam=\fivemsb

\def\dd#1{\raise1.5pt\hbox{$\,\partial\!$}/\raise-2.5pt\hbox{$\!\partial#1\,$}}

\def\tilde{\widetilde}
\def\hat{\widehat}
\def\5#1{{\mathcal #1}}

\def\RR{{\mathbb R}}
\def\CC{{\mathbb C}}

\def\NN{{\mathbb N}}

\def\ra{\rightarrow}

\def\GL{\mathop{\rm GL}\nolimits}

\def\Aut{\mathop{\rm Aut}\nolimits}

 \def\HollowBoxx #1#2#3{{\dimen0=#1 \advance\dimen0 by -#2
       \dimen1=#1 \advance\dimen1 by #3
        \vrule height 0pt depth #3 width #2
       \hskip -#3
       \vrule height #1 depth #3 width #3}}
 \def\LeftContraction{\mathord{\kern1.45pt \HollowBoxx{6pt}{3.5pt}{.4pt}}\,}

 \def\HollowBox #1#2#3{{\dimen0=#1 \advance\dimen0 by -#3
       \dimen1=#1 \advance\dimen1 by #3
        \vrule height #1 depth #3 width #3
        \vrule height 0pt depth #3 width #2
        \hskip -#3}}
 \def\RightContraction{\mathord{\, \HollowBox{6pt}{3.1pt}{.4pt}} \kern1.6pt}

\def\qed{{\hfill $\Box$}}
\newtheorem{theorem}{THEOREM}[section]

\newtheorem{lemma}[theorem]{Lemma}

\newtheorem{remark}[theorem]{Remark}
\newtheorem{proposition}[theorem]{Proposition}

\begin{document}

\begin{center}
{\Large \bf On the Kobayashi Hyperbolicity\\
\vspace{0.2cm}
of Certain Tube Domains}\footnote{{\bf Mathematics Subject Classification:} 32Q45, 31A15}
\vspace{0.4cm}\\
\normalsize Alan Huckleberry and Alexander Isaev
\end{center}

\begin{quotation} 
{\small \sl \noindent In article {\rm \cite{I2}} the second author introduced three families of tube domains in $\CC^2$ with holomorphic automorphism group isomorphic to $\RR\ltimes\RR^2$ and envelope of holomorphy equal to $\CC^2$. In the present paper we show that every domain in each of these families is Kobayashi-hyperbolic.}
\end{quotation}

\thispagestyle{empty}

\pagestyle{myheadings}
\markboth{A. Huckleberry and A. Isaev}{Hyperbolicity of Certain Tube Domains}

\setcounter{section}{0}

\section{Introduction}\label{intro}
\setcounter{equation}{0}

A connected complex manifold $X$ is called {\it Kobayashi-hyperbolic}\, if the Koba\-ya\-shi pseudodistance on $X$ is in fact a distance (see \cite{K} for details). If $X$ is equipped with a Riemannian metric, the hyperbolicity property can be stated as follows: for any point $x\in X$ there exist a neighborhood $U$ of $x$ and a constant $M>0$ such that for all holomorphic maps $f:\Delta\ra X$ with $f(0)\in U$ one has $||df(0)||<M$, where $\Delta$ is the unit disk in $\CC$ (see, e.g. \cite{L}). Verification of hyperbolicity for a particular manifold may be a difficult task. In this paper we show that certain explicitly given tube domains in $\CC^2$ are hyperbolic.

Recall that a {\it tube domain}\, in $\CC^n$ is a domain of the form $T_D:=D+i\RR^n$, where $D$ is a domain in $\RR^n$ called the {\it base}\, of $T_D$. By Bochner's theorem, the envelope of holomorphy of $T_D$ coincides with $T_{\hat D}$, where $\hat D$ is the convex hull of $D$ in $\RR^n$ (see, e.g. Section 21 in \cite{V}). For tube domains in $\CC^2$ with $T_{\hat D}\not=\CC^2$ a hyperbolicity criterion was given in \cite{L}. However, there is no reasonable sufficient condition for $T_D$ to be hyperbolic in the case $T_{\hat D}=\CC^2$. All tube domains considered in this paper fall into this last case. In particular, they do not admit any non-constant bounded holomorphic functions.

We will now introduce three families of domains in $\RR^2$ as follows:
$$
\begin{array}{l}
A_{\alpha,s,t}:=\left\{(x_1,x_2)\in\RR^2:\begin{array}{ll} x_2>t\,x_1^{\alpha}&\hbox{if $x_1>0$},\\
\vspace{-0.4cm}\\
x_2>0 & \hbox{if $x_1=0$},\\
\vspace{-0.4cm}\\
x_2>s(-x_1)^{\alpha}&\hbox{if $x_1<0$}
\end{array}\right\},\\
\vspace{-0.4cm}\\
\hspace{1.8cm}\hbox{$\alpha>0$, $\alpha\ne 1$, $s<0$, $t>0$},\\
\vspace{-0.1cm}\\
B_{s,t}:=\left\{(x_1,x_2)\in\RR^2: \begin{array}{ll} x_2>x_1\log(t\,x_1)&\hbox{if $x_1>0$},\\
\vspace{-0.4cm}\\
x_2>0 & \hbox{if $x_1=0$},\\
\vspace{-0.4cm}\\
x_2>x_1 \log(s(-x_1))&\hbox{if $x_1<0$}
\end{array}\right\},\\
\hspace{1.8cm}\hbox{$s>0$, $t>0$},\\
\end{array}
$$
$$
\begin{array}{l}
C_{\alpha,s,t}:=\Bigl\{(x_1,x_2)\in\RR^2:se^{\alpha\varphi}<r<te^{\alpha\varphi}\Bigr\},\\
\vspace{-0.4cm}\\
\hspace{1.8cm}\hbox{$\alpha>0$, $t>0$, $e^{-2\pi\alpha}t<s<t$},
\end{array}
$$
where $(r,\varphi)$ are the polar coordinates in $\RR^2$ with $\varphi$ varying from $-\infty$ to $\infty$.

If $X$ is hyperbolic, then the group $\Aut(X)$ of its holomorphic automorphisms can be given the structure of a (real) Lie group in the compact-open topology (for a discussion of this property in a more general setting see Section 1.1 in \cite{I1}). In \cite{I2} all hyperbolic manifolds $X$ with $\dim_{\CC}X=2$, $\dim_{\RR}\Aut(X)=3$ were classified up to biholomorphic equivalence, and the families of tube domains $\{T_{A_{\alpha,s,t}}\}$, $\{T_{B_{s,t}}\}$, $\{T_{C_{\alpha,s,t}}\}$ form part of the classification (they correspond to the domains listed in \cite{I2} under the headings {\bf (1)(c)},  {\bf (2)(b)},  {\bf (5)}, respectively). The holomorphic automorphism group of each of these tube domains is isomorphic to $\RR\ltimes_{\rho}\RR^2$ for some homomorphism $\rho:\RR\ra\GL(2,\RR)$. For all manifolds that appear in the classification other than the tube domains in the families $\{T_{A_{\alpha,s,t}}\}$, $\{T_{B_{s,t}}\}$, $\{T_{C_{\alpha,s,t}}\}$, verification of hyperbolicity is straightforward and therefore was omitted in \cite{I2}. In contrast, ascertaining the hyperbolicity of domains in $\{T_{A_{\alpha,s,t}}\}$, $\{T_{B_{s,t}}\}$, $\{T_{C_{\alpha,s,t}}\}$ is non-trivial. However, no statement confirming hyperbolicity for such domains was given in \cite{I2} either. In the present paper we address this issue by proving the following theorem.

\begin{theorem}\label{main}\sl Every domain in each of the families $\{T_{A_{\alpha,s,t}}\}$, $\{T_{B_{s,t}}\}$, $\{T_{C_{\alpha,s,t}}\}$ is hyperbolic.
\end{theorem}

\noindent In addition to supplementing the arguments of \cite{I2}, the proof of Theorem \ref{main} given in the next section is also of independent interest since it in fact applies to a much larger class of tube domains satisfying the condition $T_{\hat D}=\CC^2$ (see Remark \ref{moregeneral}).
\vspace{-0.1cm}\\

{\bf Acknowledgements.}  We are indebted to F. Nazarov, E. Poletsky and L. Kovalev for communicating to us an idea that has turned out to be instrumental for the proof of Theorem \ref{main} (see Remark \ref{idea} for details). This work was initiated during the first author's visit to the Australian National University in September 2011. We gratefully acknowledge support of the Australian Research Council.

\section{Proof of Theorem \ref{main}}\label{sect1}
\setcounter{equation}{0}

As pointed out in \cite{L}, for a tube domain $T_D\subset\CC^n$ the hyperbolicity property is equivalent to the following condition: for any point $x\in D$ there exist a neighborhood $U$ of $x$ in $D$ and a constant $M>0$ such that for all harmonic maps $f:\Delta\ra D$ with $f(0)\in U$ one has $||df(0)||<M$. Hence $T_D$ is not hyperbolic if and only if there exist a point $a\in D$ and a sequence $\{f_k\}$ of harmonic maps from $\Delta$ into $D$ such that $f_k(0)\ra a$ and $||df_k(0)||\ra\infty$ as $k\ra\infty$. 

Let now $D$ be a domain in one of the families $\{A_{\alpha,s,t}\}$, $\{B_{s,t}\}$, $\{C_{\alpha,s,t}\}$. Assuming that $T_D$ is not hyperbolic, we obtain a point $a=(a_1,a_2)\in D$ and a sequence $\{f_k\}$ as above, with $f_k=(u_k,v_k)$, where $u_k$, $v_k$ are real-valued harmonic functions on $\Delta$. In our proof of the theorem we utilize level sets of $u_k$. Some fundamental properties of such sets are given in the following proposition.    

\begin{proposition}\label{values} \sl For every $c\in\RR$ one has:
\vspace{-0.3cm}\\

\noindent{\rm(i)} there exists $K\in\NN$ such that $c\in u_k(\Delta)$ for all $k\ge K$,
\vspace{-0.3cm}\\

\noindent{\rm(ii)} $\hbox{\rm dist}(0,L_k(c))\ra 0$ as $k\ra\infty$, where $L_k(c):=\{z\in\Delta:u_k(z)=c\}$.  
\end{proposition}

\noindent {\bf Proof:} We first prove statement (i). Assuming it is false, we obtain a subsequence $\{u_{k_{\ell}}\}$ of the sequence $\{u_k\}$ such that for some  $c\ne a_1$ one has either $u_{k_{\ell}}<c$ in $\Delta$ (if $a_1<c$) or  $u_{k_{\ell}}>c$ in $\Delta$ (if $a_1>c$) for all $k_\ell$. Then $f_{k_{\ell}}(\Delta)$ is contained in either $D_{-}(c):=D\cap\{x_1<c\}$ or $D_{+}(c):=D\cap\{{x_1>c\}}$ for all $k_\ell$, respectively. 

Suppose first that $D$ belongs to the family $\{C_{\alpha,s,t}\}$. In this case the open sets $D_{-}(c)$ and $D_{+}(c)$ are disconnected and each of their countably many connected components is bounded. Clearly, a tube domain having a bounded base is hyperbolic. On the other hand, let $D'(c)$ be the connected component containing the point $a$. Then $f_{k_{\ell}}(\Delta)\subset D'(c)$  for large $k_\ell$, which contradicts the hyperbolicity of $T_{D'(c)}$.

Suppose next that $D$ belongs to one of the families $\{A_{\alpha,s,t}\}$, $\{B_{s,t}\}$. In this case $D_{-}(c)$ and $D_{+}(c)$ are connected. We will now show that the tube domains $T_{D_{{}_-}(c)}$ and $T_{D_{{}_+}(c)}$ are hyperbolic thus contradicting the fact that $f_{k_{\ell}}(\Delta)$ is contained in either $D_{-}(c)$ or $D_{+}(c)$ for all $k_\ell$. We use the following well-known result.

\begin{lemma}\label{eastwood}{\rm\cite{E}}\sl\, Let $X$, $Y$ be complex manifolds and $F:X\ra Y$ a holomorphic map. Suppose that $Y$ is hyperbolic and has an open cover $\{U_{\alpha}\}$ such that $F^{-1}(U_{\alpha})$ is hyperbolic for every $\alpha$. Then $X$ is hyperbolic.  
\end{lemma}

\noindent One now easily observes that $T_{D_{{}_-}(c)}$ (resp., $T_{D_{{}_+}(c)}$) is hyperbolic by choosing in Lemma \ref{eastwood} the manifold $Y$ to be $\{(z_1,0)\in\CC^2:\hbox{Re}\,z_1<c\}$ (resp., $\{(z_1,0)\in\CC^2:\hbox{Re}\,z_1>c\}$), the open cover to be $\{(z_1,0)\in\CC^2:c-m<\hbox{Re}\,z_1<c\}$ (resp., $\{(z_1,0)\in\CC^2:c<\hbox{Re}\,z_1<c+m\}$), $m\in\NN$, and the map $F$ to be the projection to the $z_1$-coordinate complex line. This completes the proof of statement (i).

We will now prove statement (ii). Assuming it is false, we obtain a subsequence $\{f_{k_{\ell}}\}$ of the sequence $\{f_k\}$ and a disk $\Delta_{r}$ of radius $0<r<1$ centered at the origin such that for some  $c\ne a_1$ the set $f_{k_{\ell}}(\Delta_r)$ is contained in either $D_{-}(c)$ (if $a_1<c$) or $D_{+}(c)$ (if $a_1>c$) for all $k_\ell$. Considering the sequence $\{\tilde f_{k_\ell}\}$ of maps from $\Delta$ to $D$ defined by $\tilde f_{k_\ell}(z):=f_{k_\ell}(rz)$ for $|z|<1$, we obtain a contradiction as in the proof of statement (i) above. The proof of Proposition \ref{values} is complete.\qed

In the remaining part of the proof of the theorem we will separately consider two cases. 

{\bf Case 1.} Suppose that $D$ belongs to one of the families $\{A_{\alpha,s,t}\}$, $\{B_{s,t}\}$. Fix $R>0$, $p>0$, $q>0$, $0<\varepsilon<1$, such that $pR>a_1$, $qR>-a_1$ and consider the open set
$$
\left\{z\in\Delta: a_1-pR<u_k(z)<a_1+qR,\,\,|z|<1-\varepsilon\right\}.
$$
For all sufficiently large $k$ the origin lies in this set, and we denote by $\Omega_k$ its connected component containing the origin. By the maximum principle for harmonic functions, $\Omega_k$ is a Jordan simply-connected domain, and we have
$$
\partial\Omega_k=\Gamma_k\sqcup\Gamma_k'\sqcup\gamma_k,
$$
where $\Gamma_k\subset L_k(a_1-pR)$, $\Gamma_k'\subset L_k(a_1+qR)$, and $\gamma_k:=\partial\Omega_k\cap\{|z|=1-\varepsilon\}$. 

For a subset $E\subset\partial\Omega_k$, let $\omega_k(E)$ be the harmonic measure of $E$ at the origin associated to $\Omega_k$. By a well-known estimate (see Theorem IV.6.2 on p. 149 in \cite{GM}) and Proposition \ref{values}, for any sufficiently large $k$ one has
$$
\omega_k(\gamma_k)\le\frac{8}{\pi}\sqrt{\frac{\hbox{\rm dist}(0,\partial\Omega_k)}{1-\varepsilon}},
$$     
which implies
\begin{equation}
\omega_k(\gamma_k)\ra 0\quad\hbox{as $k\ra\infty$}.\label{estim} 
\end{equation}
Next, let $\mu_k:=\omega_k(\Gamma_k)$ and $\mu_k':=\omega_k(\Gamma_k')$. We have
\begin{equation}
\begin{array}{l}
\displaystyle u_k(0)=\int_{\Gamma_k}u_kd\omega_k+\int_{\Gamma_k'}u_kd\omega_k+\int_{\gamma_k}u_kd\omega_k=\\
\vspace{-0.1cm}\\
\displaystyle \hspace{1.5cm}(a_1-pR)\mu_k+(a_1+qR)\mu_k'+\int_{\gamma_k}u_kd\omega_k.
\end{array}\label{integral}
\end{equation}
Since on $\gamma_k$ the function $u_k$ is bounded from above and below by constants independent of $k$,  from (\ref{estim}) we obtain that the last summand in (\ref{integral}) tends to zero as $k\ra\infty$. Thus, (\ref{estim}) and (\ref{integral}) yield
$$
\begin{array}{l}
\mu_k+\mu_k'\ra 1,\\
\vspace{-0.1cm}\\
(a_1-pR)\mu_k+(a_1+qR)\mu_k'\ra a_1
\end{array}\quad\hbox{as $k\ra\infty$},
$$
which implies
\begin{equation}
\mu_k\ra\frac{q}{p+q},\quad \mu_k'\ra \frac{p}{p+q}\quad\hbox{as $k\ra\infty$}.\label{mupq}
\end{equation}

We will now consider two situations.

{\bf Case 1.1.} Assume that $D=A_{\alpha,s,t}$ for some $\alpha>0$, $\alpha\ne 1$, $s<0$, $t>0$. Then
$$
\begin{array}{l}
\hbox{$v_k(z)>s(pR-a_1)^{\alpha}$ for $z\in\Gamma_k,\,\gamma_k$},\\
\vspace{-0.1cm}\\
\hbox{$v_k(z)>t(a_1+qR)^{\alpha}$ for $z\in\Gamma_k'$}.
\end{array}
$$
Therefore, we have
\begin{equation}
\begin{array}{l}
\displaystyle v_k(0)=\int_{\Gamma_k}v_kd\omega_k+\int_{\Gamma_k'}v_kd\omega_k+\int_{\gamma_k}v_kd\omega_k\ge\\
\vspace{-0.1cm}\\
\displaystyle\hspace{1.5cm}s(pR-a_1)^{\alpha}\mu_k+t(a_1+qR)^{\alpha}\mu_k'+s(pR-a_1)^{\alpha}\omega_k(\gamma_k).
\end{array}\label{integral2}
\end{equation}
Letting in the above inequality $k\ra\infty$ and using (\ref{estim}), (\ref{mupq}), we then obtain
\begin{equation}
\begin{array}{l}
\displaystyle a_2\ge s(pR-a_1)^{\alpha}\frac{q}{p+q}+t(a_1+qR)^{\alpha} \frac{p}{p+q}=\\
\vspace{-0.1cm}\\
\displaystyle\hspace{4cm}\frac{R^{\alpha}}{p+q}\left[sq\left(p-\frac{a_1}{R}\right)^{\alpha}+tp\left(q+\frac{a_1}{R}\right)^{\alpha}\right].
\end{array}\label{estim1}
\end{equation}
Choosing $p,q$ such that $(q/p)^{\alpha-1}>|s|/t$ and letting $R\ra\infty$, we now observe that the right-hand side of (\ref{estim1}) can be made arbitrarily large. This contradiction concludes the proof of the theorem in the case when $D$ belongs to the family $\{A_{\alpha,s,t}\}$.

{\bf Case 1.2.} Assume that $D=B_{s,t}$ for some $s>0$, $t>0$. Then
$$
\begin{array}{l}
\hbox{$v_k(z)>(a_1-pR)\log(s(pR-a_1))$ for $z\in\Gamma_k$},\\
\vspace{-0.1cm}\\
\hbox{$v_k(z)>(a_1+qR)\log(t(a_1+qR))$ for $z\in\Gamma_k'$},\\
\vspace{-0.1cm}\\
\hbox{$v_k(z)>C$ for $z\in\gamma_k$, where $C$ is a constant independent of $k$.}
\end{array}
$$ 
Therefore, analogously to (\ref{integral2}) we have
$$
\begin{array}{l}
v_k(0)\ge (a_1-pR)\log(s(pR-a_1))\mu_k+\\\vspace{-0.1cm}\\
\hspace{6cm}(a_1+qR)\log(t(a_1+qR))\mu_k'+C\omega_k(\gamma_k).
\end{array}
$$
Letting in the above inequality $k\ra\infty$ and using (\ref{estim}), (\ref{mupq}), we then obtain
\begin{equation}
\begin{array}{l}
\displaystyle a_2\ge (a_1-pR)\log(s(pR-a_1))\frac{q}{p+q}+\\
\vspace{-0.3cm}\\
\displaystyle\hspace{6cm}(a_1+qR)\log(t(a_1+qR))\frac{p}{p+q}=\\
\vspace{-0.3cm}\\
\displaystyle\hspace{1cm}\frac{R}{p+q}\left[q\left(\frac{a_1}{R}-p\right)\log(s(pR-a_1))+\right.\\
\vspace{-0.3cm}\\
\displaystyle\hspace{6cm}\left.p\left(\frac{a_1}{R}+q\right)\log(t(a_1+qR))\right].
\end{array}\label{estim2}
\end{equation}
Choosing $p,q$ such that $tq>sp$ and letting $R\ra\infty$, we now observe that the right-hand side of (\ref{estim2}) can be made arbitrarily large. This contradiction concludes the proof of the theorem in the case when $D$ belongs to the family $\{B_{s,t}$\}.

{\bf Case 2.} Suppose now that $D$ belongs to the family $\{C_{\alpha,s,t}\}$. Fix $c>a_1$, $0<\varepsilon<1$ and consider the open set
$$
\left\{z\in\Delta: u_k(z)<c,\,\,|z|<1-\varepsilon\right\}.
$$
For all sufficiently large $k$ the origin lies in this set, and we denote by $\Omega_k$ its connected component containing the origin. As in Case 1, $\Omega_k$ is a Jordan simply-connected domain, and we have
$$
\partial\Omega_k=\Gamma_k\sqcup\gamma_k,
$$
where $\Gamma_k\subset L_k(c)$ and $\gamma_k:=\partial\Omega_k\cap\{|z|=1-\varepsilon\}$. 

Recall from the proof of Proposition \ref{values} that the open set $D_{-}(c)=D\cap\{x_1<c\}$ has countably many connected components and each component is bounded. Let $D'(c)$ be the connected component of $D_{-}(c)$ containing $a$. Clearly, $f_k(\Omega_k)\subset D'(c)$ for all sufficiently large $k$. This implies that on $\gamma_k$ the function $u_k$ is bounded from below by a constant independent of $k$ if $k$ is sufficiently large. 

For a subset $E\subset\partial\Omega_k$, we let $\omega_k(E)$ be the harmonic measure of $E$ at the origin associated to $\Omega_k$ and $\mu_k:=\omega_k(\Gamma_k)$. Arguing as in Case 1, we then see that (\ref{estim}) holds, that is, $\mu_k\ra 1$ as $k\ra\infty$. 

Next, we have
\begin{equation}
\displaystyle u_k(0)=\int_{\Gamma_k}u_kd\omega_k+\int_{\gamma_k}u_kd\omega_k=
c\mu_k+\int_{\gamma_k}u_kd\omega_k.\label{integral3}
\end{equation}
Since on $\gamma_k$ the function $u_k$ is bounded from above and below by constants independent of $k$,  from (\ref{estim}) we obtain that the last summand in (\ref{integral3}) tends to zero as $k\ra\infty$. Thus, (\ref{integral3}) implies $a_1=c$, which contradicts our choice of $c$. This completes the proof of the theorem.\qed

\begin{remark}\label{moregeneral} \rm The proof of Theorem \ref{main} in fact applies to more general domains. Indeed, let $D$ be a domain of the form $\{(x_1,x_2):x_2>h(x_1)\}$, where $h\in C(\RR)$ and satisfies the following property: for every $b\in\RR$ there exist $p>0$, $q>0$ such that
$$
qh(b-pR)+ph(b+qR)\ra \infty\quad\hbox{as $R\ra\infty$}.
$$
Then the argument given for Case 1 yields that $T_D$ is hyperbolic. Next, let $D$ be a domain bounded by two general spirals, where a spiral is a curve defined by the equation $r=g(\varphi)$, with $g$ being an increasing function of $\varphi$ such that $\lim_{\varphi\ra-\infty}g(\varphi)=0$ and $\lim_{\varphi\ra+\infty}g(\varphi)=\infty$. Then the argument given for Case 2 shows that $T_D$ is hyperbolic.    
\end{remark}

\begin{remark}\label{idea} \rm Before attempting to prove Theorem \ref{main} in full generality, we set out to show that the domain $T_{A_{{}_{3,-1,1}}}$ is Brody hyperbolic (recall that $A_{3,-1,1}=\{(x_1,x_2)\in\RR^2: x_2>x_1^3\}$). Brody hyperbolicity for a tube domain is equivalent to the non-existence of a non-constant harmonic map $f=(u,v)$ from the plane into the base of the domain (cf. \cite{L}). Regarding this question, F. Nazarov suggested that we consider the connected component $\Omega(R,\rho)$ containing the origin of the open set
$$
\left\{z\in\CC: -R<u(z)<2R,\,\,|z|<\rho\right\}
$$
(assuming without loss of generality that $u(0)=0$), where $R>0$ and $\rho>0$ are large. Then the harmonic measure at the origin of the portion of $\partial\Omega(R,\rho)$ where $|z|=\rho$ tends to 0 as $\rho\ra\infty$, and letting $R\ra\infty$ one obtains that $v(0)$ is estimated from below by an arbitrarily large number. With F. Nazarov's kind permission, we used a similar approach in Case 1 of our proof of Theorem \ref{main}. We also point out that the idea to consider harmonic measures associated to domains bounded by level sets of $u$ was independently suggested to us by E. Poletsky and L. Kovalev.   
\end{remark}

{\obeylines

\noindent Alan Huckleberry:
\noindent Fakult\"at f\"ur Mathematik
\noindent Ruhr-Universit\"at Bochum
\noindent Universitatstra\ss e 150
\noindent 44801 Bochum
\noindent Germany
\noindent and
\noindent School of Engineering and Science
\noindent Jacobs University Bremen
\noindent Campus Ring 1
\noindent 28759 Bremen
\noindent Germany
\noindent e-mail: ahuck@gmx.de
\vspace{0.3cm}

\noindent Alexander Isaev:
\noindent Department of Mathematics
\noindent The Australian National University
\noindent Canberra, ACT 0200
\noindent Australia
\noindent e-mail: alexander.isaev@anu.edu.au
}

\end{document}